\documentclass{article}
\usepackage{subfigure}
\usepackage{graphicx,amssymb,amsthm,amsmath}

\title{ A Graph-Theoretic Approach to Efficient Voice-Leading}

\newtheorem{definition}{Definition}

\author{Susannah Wixey and Rob Sturman$^{\ast}$\thanks{$^\ast$Corresponding author. Email: r.sturman@leeds.ac.uk
\vspace{6pt}}\\\vspace{6pt} {\em{School of Mathematics, University of Leeds, Leeds LS2 9JT, UK}} }

\begin{document}

\maketitle

\begin{abstract}
We study the Neo-Riemannian principle of parsimonious voice-leading using tools and techniques from classical graph theory and the modern field of complex networks. We quantify the relative importance of particular chords within this framework. The graph-theoretic notion of eccentricity suggests that when working in a harmonic scheme dictated by any common musical scale, no triad is any more isolated than any other. Complex network theory refines this idea, and in this context provides measures of how important particular triads might be for the  flow of chord progressions through the harmonic network. We review and compare several different such measures of centrality and communicability.
\end{abstract}

\section{Introduction}

Neo-Riemannian theory is based on the idea of connecting chords according to some definition of harmonic proximity. This arose from a reappraisal of the work of Hugo Riemann by music theorists including Richard Cohn, David Lewin, Dmitri Tymoczko and others~\cite{cohn1998introduction}. The premise of harmonic proximity requires a notion of a metric, or distance, between chords, making an appeal to mathematical formalism natural. There is a variety of different ways to define such chordal distance, as discussed in \cite{cohn2012audacious}. A fundamental idea is that of {\em parsimonious voice-leading}. This essentially describes an efficiency of movement between chords, illustrated by the simple example of the proximity of $C$ major and $e$ minor triads. Both share a pair of common pitches ($E$ and $G$), while their remaining pitches are neighbours in the chromatic scale, being only a semitone apart. Choosing an initial scale, or set of pitches, thus defines in turn a set of available triads which exist in the set of pitches, and a set of parsimonious voice-leading connections. These triads and connections, viewed as vertices and edges, naturally lead to a mathematical treatment of the network produced. 

One successful view of this type of structure is the geometrical. Many pitch sets, including the diatonic and the chromatic scales, produce a well-defined geometry, which lends itself well to musical analysis. This view has been comprehensively studied in \cite{tymoczko2011geometry}, which  explores geometric features and how they constrain or inspire composers who wish to create music with a sense of tonality. A drawback of the geometrical approach is that a given pitch set might produce a network with no accessible geometry, for example, the network emanating from a harmonic minor scale. Interpreting the network as a mathematical graph rather than a geometry allows a greater range of pitch sets to be studied. A graph theoretic approach has been taken in work such as \cite{douthett1998parsimonious} and \cite{waller1978some}, in which the focus is on pitch sets which produce networks possessing some symmetry. In the present paper we do not look for any geometrical or symmetrical structure in the networks, but ask purely graph theoretic questions, such as how to quantify connectivity and traversability. 

Classical graph theory provides some such measures, but recent work in complex network theory asks different questions. In particular, how efficiently can information flow through a network? Typically a network in complexity studies is huge --- far, far bigger than the networks we study here --- but still the principles and diagnostics, such as {\em centrality} and {\em communicability}, developed for such large graphs have a pertinent interpretation to voice-leading graphs, providing information about how important particular chords are in optimising the efficiency of a chord progression between a pair of triads in a given harmonic landscape.

\subsection{Notation} 

Throughout this paper we use the following notation: major triads are denoted with upper-case script letter; minor triads in lower case. Augmented triads are given in upper-case with a superscript $+$; diminished triads are lower-case with a superscript $o$. Thus, $F\sharp$, $f\sharp$, $F\sharp^+$, $f\sharp^o$ denote major, minor, augmented and diminished triads respectively, on the pitch $F\sharp$.

\section{Voice Leading Graphs}

Parsimonious, or efficient, voice-leading is a notion of counterpoint which is prevalent in Western tonal music and which is frequently discussed in relation to neo-Riemannian theory. Following \cite{tymoczko2011geometry}, we begin with some chosen scale (or more generally, a collection of pitches) and determine all possible major, minor, augmented and diminished triads which can be constructed using the notes of our collection. For the purposes of this paper we will consider the different inversions of a triad to be equivalent. We will also accept octave equivalence, so we essentially consider pitch classes rather than pitches and unless otherwise stated we will be working in twelve-tone equal temperament. Next, we construct a mathematical graph, with each vertex representing a triad and two vertices being joined by an edge if the corresponding triads can be connected by `efficient' voice leading. In this paper we consider `single-step voice-leading' where we require that two notes remain unchanged and the third note may move by a single step {\em within the scale or collection of notes} in which we are working, however large that step may be. That is, two vertices of our graph are connected precisely when the corresponding triads contain two common tones and the remaining notes are a single `scale-step' apart. This model is restrictive, but well-defined. It serves as a basic procedure to illustrate principles which could then be generalised to other situations: for example we might: allow chords of more than three notes; allow more than one note to move at a time,  or for one note to move by several steps; work in $n$-tone equal temperament, or other tuning systems. In all cases the graphs obtained by this method will be simple (that is, they contain no loops or multiple edges), and the graphs are undirected. Some generalisations might call for a weighted graph (that is, with a numerical value assigned to each edge), or a directed graph (that is, with directions specified on edges, so that a chord progression may be permitted in one direction, but not the reverse), or a multigraph (that is, with multiple edges permitted between a pair of vertices), and the associated theory.

\subsection{Construction}

Here we give an explicit set-theoretic description of the procedure of constructing  voice-leading graph. Consider the (ordered) pitch class set $\mathcal{C} = \{ 0 , 1 , 2 , \ldots , 11\}$, representing $12$-tone equal temperament, with $0$ representing the pitch class $C$, $1$ representing $C\sharp$, and so on. Choose an ordered set $\mathcal{P} = \{ p_1, p_2, \ldots , p_m \}$ of $m$ pitch classes, with $p_i \in \mathcal{C}$. There are ${12 \choose m}$ possible choices for $\mathcal{P}$. A major triad on the pitch $p_j$ exists in $\mathcal{P}$ if the set $\{ p_j ,  p_j+4, p_j+7\} \subseteq \mathcal{P}$, with addition taken modulo $12$. Similarly, minor, diminished and augmented triads exist if $\{ p_j ,  p_j+3, p_j+7\} \subseteq \mathcal{P}$, $\{ p_j ,  p_j+3, p_j+6\} \subseteq \mathcal{P}$ and $\{ p_j ,  p_j+4, p_j+8\} \subseteq \mathcal{P}$ respectively. We label the collection of all such triads $\mathcal{V} = \mathcal{V}(\mathcal{P}) = \{ v_1 , v_2, \ldots , v_n\}$. These will form the $n$ vertices of the voice-leading graph.

An edge exists between vertices $v_r = \{ p_{r_1},p_{r_2},p_{r_3} \}$ and $v_s = \{ p_{s_1},p_{s_2},p_{s_3} \}$ if and only if $v_r \cap v_s$ contains exactly two elements, and iff the two elements of the symmetric difference $v_r \triangle v_s$ are neighbours (modulo $n$) in $\mathcal{P}$, that is, $\nexists \hat{p} \in \mathcal{P} \mbox{ such that } p'<\hat{p}<p'' \mbox{ for } p',p'' \in v_r \triangle v_s$.

For example, let $\mathcal{P} = \{ 0,2,4,5,7,9,11\}$, which contains the major triads $C = \{ 0,4,7\}$, $F = \{ 5,9,0\}$  and $G = \{ 7,11,2\}$, the minor triads $d = \{ 2,5,9\}$, $e = \{ 4,7,11\}$ and $a = \{ 9,0,4\}$, and the diminished triad $b^o = \{ 11,2,5\}$. There are no augmented triads. The triad pair $(C,e)$ are connected since $C \cap e = \{4,7\}$ and $C\triangle e = \{ 0,11\}$ contains neighbours in the ordered set $\mathcal{P}$. Similarly, triad pairs $(e,G)$, $(G,b^o)$, $(b^o,d)$, $(d,F)$, $(F,a)$ and $(a,C)$ are connected, but no other triad pairs.

\begin{figure}
\centering
\subfigure[Diatonic scale, $\{0,2,4,5,7,9,11\}$.]{
\includegraphics[width=0.45\linewidth]{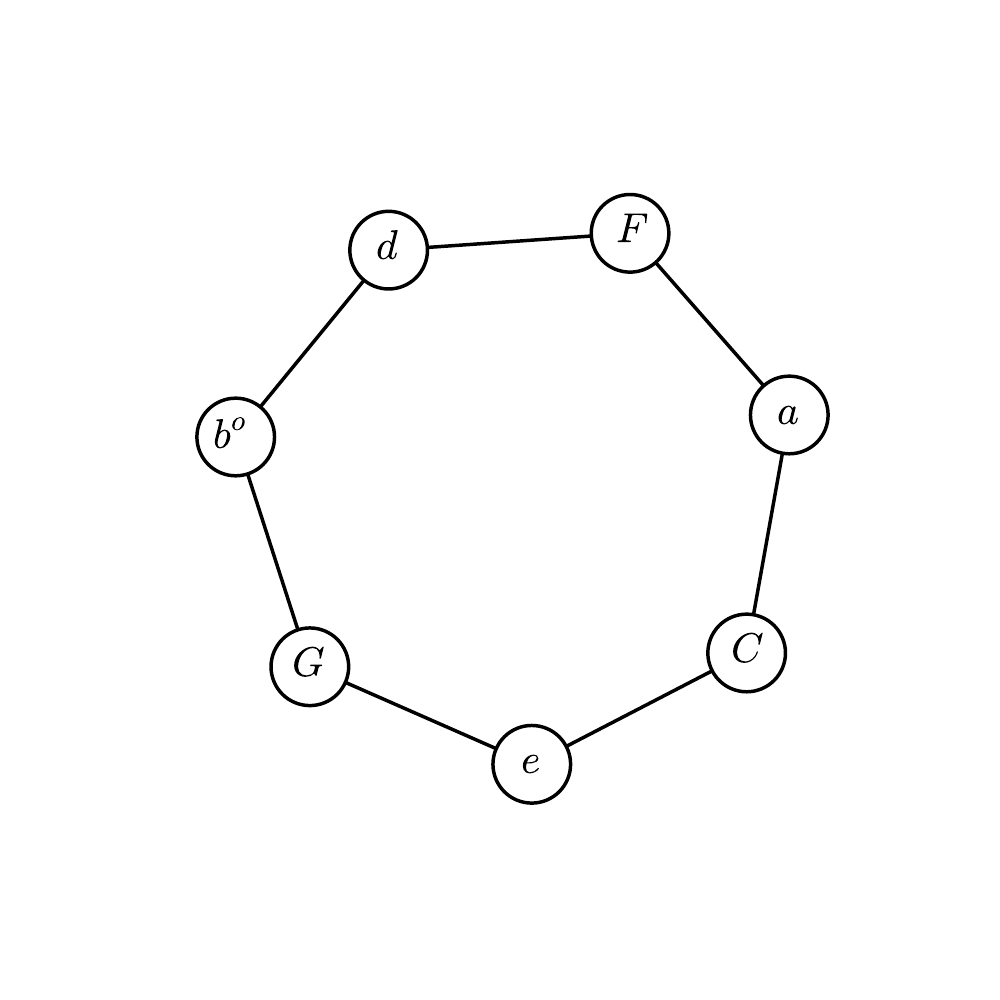}\label{fig:diatonic}}
\subfigure[Hexatonic scale, $\{ 0,1,4,5,8,9\}$.]{
\includegraphics[width=0.45\linewidth]{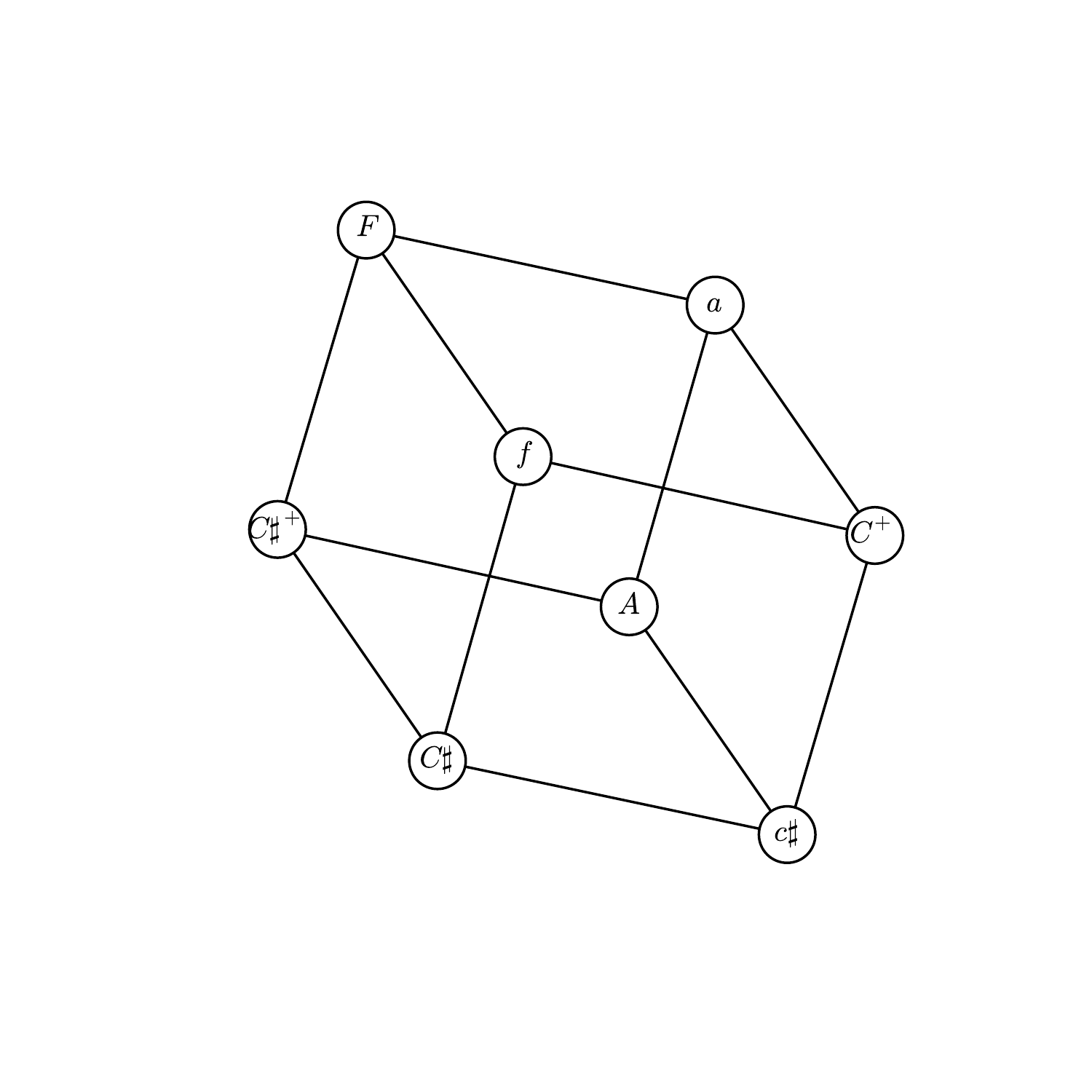}\label{fig:hexatonic}}
\caption{Voice-leading relations among triads contained in particular musical scales, represented as mathematical graphs. Vertices represent triads, and edges represents parsimonious voice-leading.}\label{fig:diatonic_graph}
\end{figure}

The resulting graph, shown in figure \ref{fig:diatonic}, is naturally identical to the corresponding figure in texts such as \cite{tymoczko2011geometry}, but here we focus on graph-theoretic  rather than geometrical properties. Of course, when $\mathcal{P}$ consists of the notes of any major (or natural minor) scale, a similar circular graph is produced. Other choices for $\mathcal{P}$ may give different graphs. Figure \ref{fig:hexatonic} shows the voice-leading graph for the notes of the hexatonic scale $\mathcal{P} = \{ 0, 1,4,5,8,9\}$, which produces a graph with eight vertices, each connected to two others. Note that this graph has a different structure to the {\em HexaCycles} of \cite{douthett1998parsimonious} as we are allowing augmented triads. Geometrically, figure \ref{fig:hexatonic} depicts a cube, but again we concentrate on graph-theoretic properties.



\begin{figure}
\centering
\includegraphics[width=0.8\linewidth]{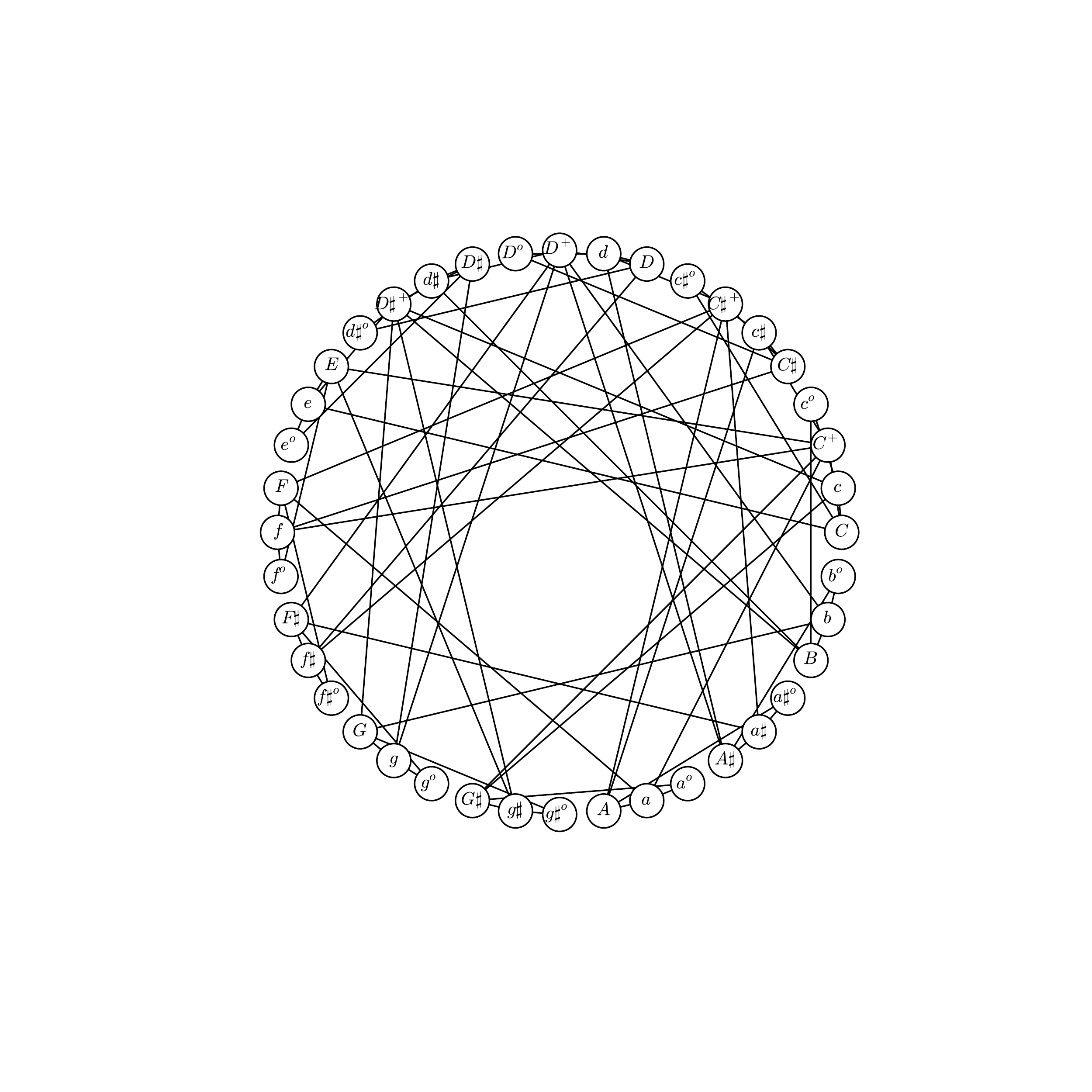}
\caption{Voice-leading relations among all forty major, minor, augmented and diminished triads available in the chromatic scale in 12-tone equal temperament. This figure is similar to a flattened version of the {\em Douthett-Steinbach Chicken Wire Torus} (figure 7 of \cite{douthett1998parsimonious}), previously discovered by \cite{waller1978some}, although these do not include augmented or diminished triads. }\label{fig:chrom_graph}
\end{figure}

\begin{figure}
\centering
\subfigure[$\mathcal{P} =\{0,2,4,5,6,7,9,11\}$]{
\includegraphics[width=0.45\linewidth]{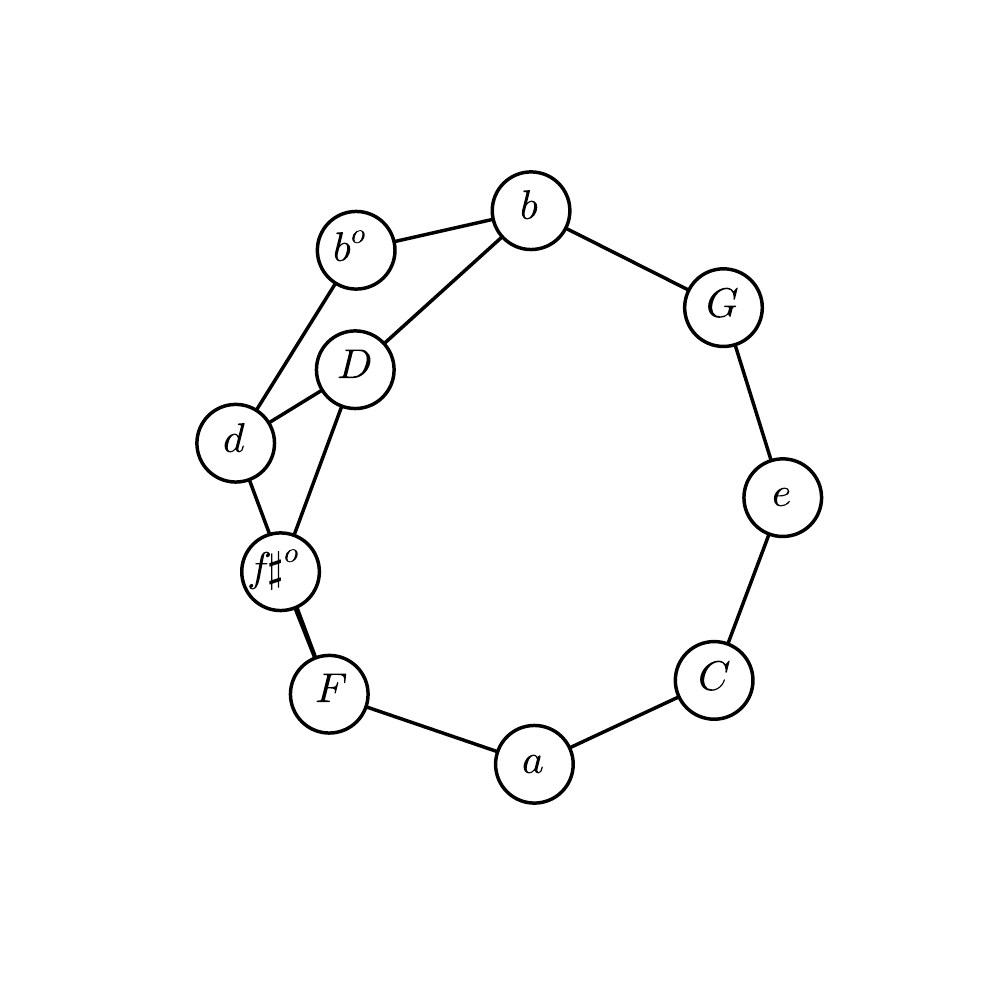}\label{fig:Cmaj_Fsharp}}
\subfigure[$\mathcal{P} = \{0,2,4,5,7,8,9,11\}$]{
\includegraphics[width=0.45\linewidth]{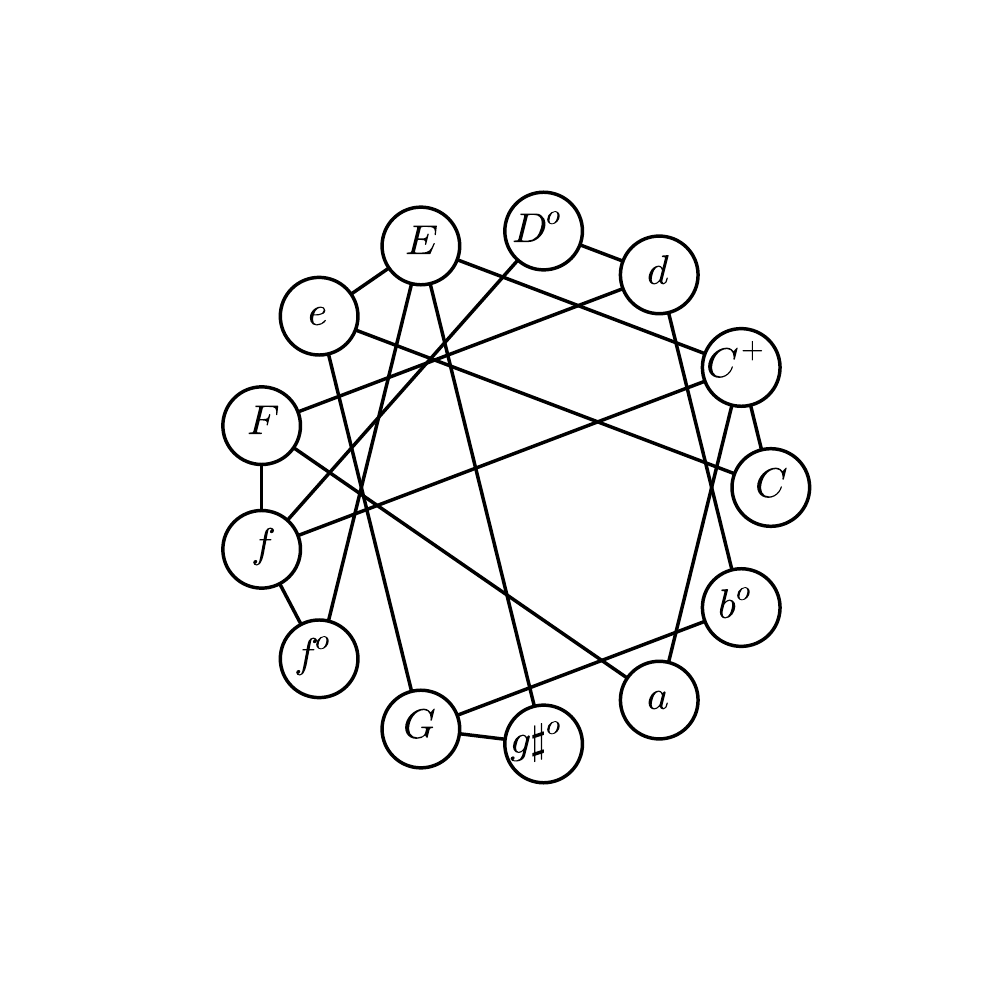}\label{fig:Cmaj_Gsharp}}
\caption{Voice-leading relations among triads constructable from a diatonic scale on $C$, plus a single extra black note.}\label{fig:Cmaj_plus}
\end{figure}

When $\mathcal{P} = \mathcal{C}$, and we can form all triads on all notes in 12-tet, the difference between the geometric and graph-theoretic viewpoints is more apparent. Figure \ref{fig:chrom_graph} show the voice-leading graph for the chromatic scale. In \cite{tymoczko2011geometry} this is interpreted as a twisted two-torus in three dimensions, whereas here we sacrifice the geometry in favour of accessing graph-theoretic diagnostics which quantify how well-connected, and how easily traversed, is the underlying harmonic space. We note that similar graphs have been previously studied
This model also lends itself well to those scales with various numbers of notes for which there is not always an obvious geometrical interpretation. For example, figure \ref{fig:Cmaj_plus} shows the pair of graphs produced when the pitch set $\mathcal{P}$ is the $C$ major scale extended first with the note $F\sharp$, and second with the note $G\sharp$. Each pitch set contains exactly 8 pitches, but different numbers of triads and efficient progressions are possible. Accepting these pitchs sets as possible harmonic palettes, it is reasonable to wish to quantify the differences between the resulting graphs in some meaningful way.


\section{Classical Graph Theory}

To analyse these voice-leading graphs we utilise common graph-theoretic properties. These make quantifiable distinctions between different graphs, and, musically,  provide information about the level of harmonic choice available to composers wishing to employ these contrapuntal techniques. We note that some graph-theoretic properties, such as connectivity, colouring and group structure have been studied by \cite{walton2010graph}, but diminished and augmented triads are not considered in that paper. 

\subsection{Vertices and  Edges}

A mathematical graph $G$ consists of a vertex set $V$ together with a set $E$ containing edges which may link some vertices.  Since some of the terminology in this area varies between authors, we begin with some basic definitions. Where there is no risk of ambiguity we will use juxtaposition notation, denoting an edge connecting two vertices, $u$ and $v$, by $uv$.

\begin{definition}
The \textbf{degree} of a vertex is the number of edges with that vertex as an endpoint, and a graph is said to be \textbf{regular} if all vertices have the same degree.
\end{definition}

Musically, if a voice-leading graph is regular, this tells us that the level of harmonic choice is consistent throughout. For example, the hexatonic scale known sometimes as the `symmetrical augmented' scale, which is constructed by alternating between intervals of a minor third and a semitone, produces a regular graph, shown in figure~\ref{fig:hexatonic}. Similarly, the diatonic scale graph of figure~\ref{fig:diatonic} is regular, but the chromatic graph in figure~\ref{fig:chrom_graph} and those in figure~\ref{fig:Cmaj_plus} are not.

\begin{definition} A \textbf{walk} or \textbf{edge sequence} is a sequence of vertices $v_1,v_2,...,v_{k-1},v_k$ such that $v_iv_{i+1}$ is an edge for all $1 \leq i \leq k-1$.
\end{definition}

Thus a walk is simply a chord progression between $k$ triads. 

\begin{definition} A graph with vertex set $V=\{ v_1,v_2,...,v_n \}$ can be represented by an $n \times n$ \textbf{adjacency matrix}, $A=(a_{ij})$, the entries of which are given by 
$$
 a_{ij}=\begin{cases}
    1 & \text{if there is an edge joining $v_i$ and $v_j$},\\
    0 & \text{otherwise}.
  \end{cases}
$$
\end{definition}

The $ij$-th entry of the $k$th power of adjacency matrix $(A^k)_{ij}$ represents the number of walks of length $k$ from $v_i$ to $v_j$. This gives us the number of different chord progressions of length $k$ between specific pairs of triads.

\begin{definition} Let $G$ and $G'$ be graphs with vertex sets $V$ and $V'$ respectively. $G$ and $G'$ are \textbf{isomorphic} if there exists a bijective mapping $\theta:V\rightarrow V'$ such that for all $u$,$v\in V$, $u$ and $v$ are joined by an edge in $G$ if and only if $\theta(u)$ and $\theta(v)$ are joined by an edge in $G'$. \end{definition}

If two graphs are isomorphic they have the same number of vertices and edges and the same degree distribution.

Graphs which are isomorphic exhibit the same structure, so when two voice-leading graphs are isomorphic this means we have the same number of available triads, and the same amount of harmonic choice available to the composer. However, two graphs may be isomorphic and contain different \textit{types} of chords as we see in the following example. The four graphs obtained by starting with a $C$ major scale and adding, respectively, the pitches $C\sharp, D\sharp, F\sharp$ or $A\sharp$ are all isomorphic graphs on ten vertices (one is shown in figure~\ref{fig:Cmaj_Fsharp}). However, we do not have the same types of chords: the pitch sets $\{ 0,1,2,4,5,7,9,11\}$ and $\{0,2,3,4,5,7,9,11\}$ each contain an augmented triad whereas the pitch sets $\{0,2,4,5,6,7,9,11\}$ and $\{0,2,4,5,7,9,10,11\}$ do not. The graph obtained by adding the pitch $G\sharp$ to the $C$ major scale has thirteen vertices (see figure~\ref{fig:Cmaj_Gsharp}) and hence is not isomorphic to the others. In fact, the four smaller graphs are isomorphic to a subgraph of the $\{0,2,4,5,7,8,9,11\}$ network, revealing the greater harmonic choice in this scale.


\subsection{Distance and traversability}

Since we are primarily concerned with exploring efficiency of voice leading, it is clearly useful to consider how one might traverse these graphs.

\begin{definition} A \textbf{path} is a walk in which no vertex is repeated. A walk in which the initial and final vertices coincide (that is, $v_1=v_k$), but no other vertex is repeated,  is called a \textbf{circuit}.
\end{definition}

Walks, paths and circuits all represent chord progressions through the graphs. Paths allow each chord to be used only once within the progression, whereas walks bear no such restrictions. Circuits require the progression to begin and end on the same chord  (so one might wish to consider circuits beginning on the tonic triad).

\begin{definition}The \textbf{(geodesic) distance} between two vertices, $v_i$ and $v_j$, is defined to be the number of edges in the shortest path which connects them. We denote this by $d(v_i,v_j).$
\end{definition}

The distance between a pair of vertices can be thought of as the `size' of the voice-leading between two triads since it measures how many steps must be moved in total between the chords.  There are various metrics used for measuring voice-leading size \cite{tymoczko2011geometry,cohn2012audacious}, arguably none of them standard. Thus, we adopt the graph theoretic notion of distance as our definition of voice-leading size, remaining aware of the limitations of this definition. For example, we cannot recognise so-called `voice-crossings' which intuitively ought to make the voice-leading size larger. One might also wish to consider, not only path length, but chromatic length, that is, how many semitones have moved, rather than scale-steps. This could be achieved by weighting each edge to show how much movement has occurred in terms of semitones. This would then allow for the same graph theoretic analysis, but using the analogous tools designed for weighted graphs.

\begin{definition}
The \textbf{eccentricity} $e(v)$ of a vertex $v$ is the greatest geodesic distance of $v$ from any other vertex.
\end{definition}

\begin{definition} 
The \textbf{diameter} of $G$, $d(G)$ is the greatest distance between any two vertices of $G$. That is, $d(G) = \max_{v\in V} e(v)$. The \textbf{radius} of $G$ is the minimum eccentricity of any vertex, that is, $r(G) = \min_{v \in V} e(v)$. 
\end{definition}
The diameter can be thought of as representing the maximum amount of movement that must occur to create a voice-leading between any two triads, in other words, the longest necessary chord progression. 

\begin{definition}
A \textbf{central vertex} is any vertex  $v$ which achieves the radius, that is $e(v)=r$.  A \textbf{peripheral} vertex is any vertex $v$ which achieves the diameter, that is, $e(v)=d$. A graph is \textbf{self-centred} if $d(G)=r(G)$, that is, if every vertex has the same eccentricity.
\end{definition}

In a self-centred graph every vertex is both central and peripheral, so no vertex is any more isolated than any other. There are $\sum_{k=3}^{12}  {12 \choose k } = 4017$ possible pitch sets of between three (so there are sufficient pitches for at least one triad) and twelve pitches. Of these, 642 produce empty graphs (that is, there are no major, minor, augmented or diminished triads contained in $\mathcal{P}$), and 2 (the whole tone scales $\mathcal{P} = \{ 0,2,4,6,8,10\}$ and $\mathcal{P}  \{ 1,3,5,7,9,11\}$) produce unconnected graphs, with two vertices and no edges. Of the remaining 3373, 1857 produce self-centred graphs, while 1516 produce non-self-centred graphs.

Notably, standard musical scales, including major (and therefore natural minor), harmonic minor, melodic minor, octatonic, hexatonic and chromatic scales, produce self-centred graphs. So indeed do all the previously discussed scales produced from adding any combination of black notes to the $C$ major scale.  The musical implication is that in such harmonic frameworks, no chord is more or less isolated than any other: each triad requires the same length chord progression to reach its most distant triad. The list of non-self-centred scales include many with large gaps between elements of $\mathcal{P}$, but there exist some (less standard) musical scales that are non-self-centred. For example, the Mixolydian Augmented scale $\{0,2,4,5,8,9,10\}$ has a diameter of 4 and a radius of 3, while the Enigmatic minor scale (reputedly invented by Guiseppe Verdi) $\{ 0,1,4,6,8,10,11\}$ forms a line graph of five triads, and therefore is not self-centred.

Geodesic distances give a crude first impression of how difficult a graph is to traverse. These can be refined using some very well-known concepts from graph theory. We next discuss Hamiltonian paths and circuits, which relate closely to one of the most famous and computationally expensive extant problems of computer science: the Travelling Salesperson Problem. We then examine Eulerian trails, discussed by the eponymous Euler in his solution to the Bridges of K\"{o}nigsberg problem, the solution of which is arguably the oldest theorem of graph theory.

\begin{definition} A \textbf{Hamiltonian circuit} is a circuit which visits every vertex of a graph.\end{definition}

A Hamiltonian circuit represents a voice-leading in which every available triad is used precisely once. Although this is clearly very restrictive, the principle could potentially be used by a composer wishing to exhibit all triads as efficiently as possible. Consider once more the graph of the hexatonic scale shown in figure \ref{fig:hexatonic}.  Due to the relationship between this graph and the cube, its properties have been explored in depth --- it contains no less than 12 Hamiltonian circuits, that is, there are 12 distinct ways to write a chord progression obeying parsimonious voice-leading, which start and end on the same triad and includes each other exactly once.  In general it is nontrivial to discover whether a graph contains a Hamiltonian path or circuit, since these are both NP-complete problems \cite{karp1972reducibility}. Hamiltonian cycles in the context of chordal and timbral musical morphologies have recently been studied in \cite{akhmedov2014chordal}.



The Hamiltonian property ensures all triads are visited, but does not require every edge to be used. A composer might wish to use every possible progression between triad pairs as efficiently as possible. This is described by the Eulerian property.

\begin{definition} An \textbf{Eulerian trail} begins and ends with the same vertex and uses each edge of a graph precisely once, but may repeat vertices. A \textbf{semi-Eulerian trail} is like an Eulerian trail except that it begins and ends on different vertices. A graph is called \textbf{Eulerian} if it contains an Eulerian trail and \textbf{semi-Eulerian} if it contains a semi-Eulerian trail but not an Eulerian trail.
\end{definition}

\begin{figure}
\subfigure[]{
\includegraphics[width=0.45\linewidth]{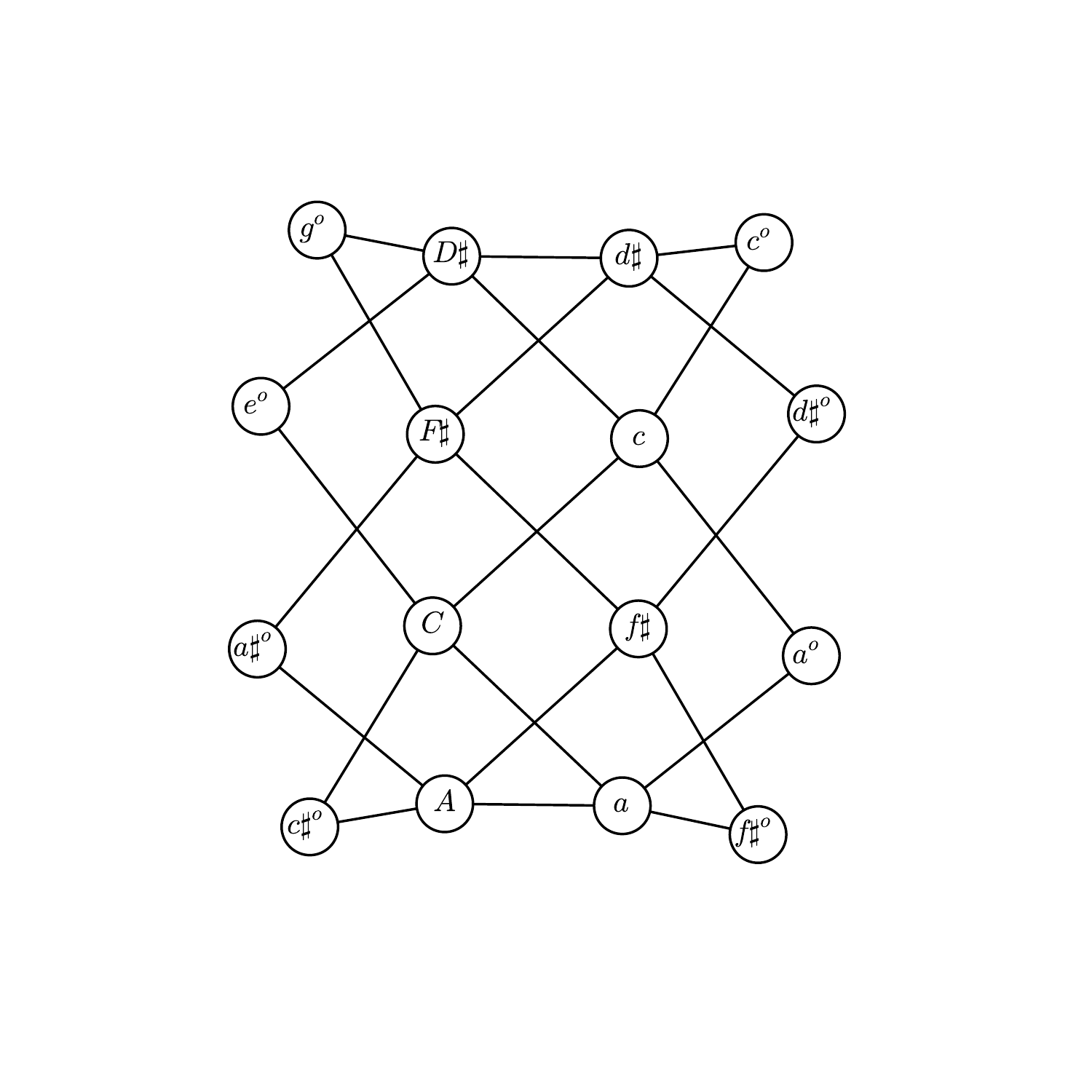}\label{fig:octatonic1}}
\subfigure[]{
\includegraphics[width=0.45\linewidth]{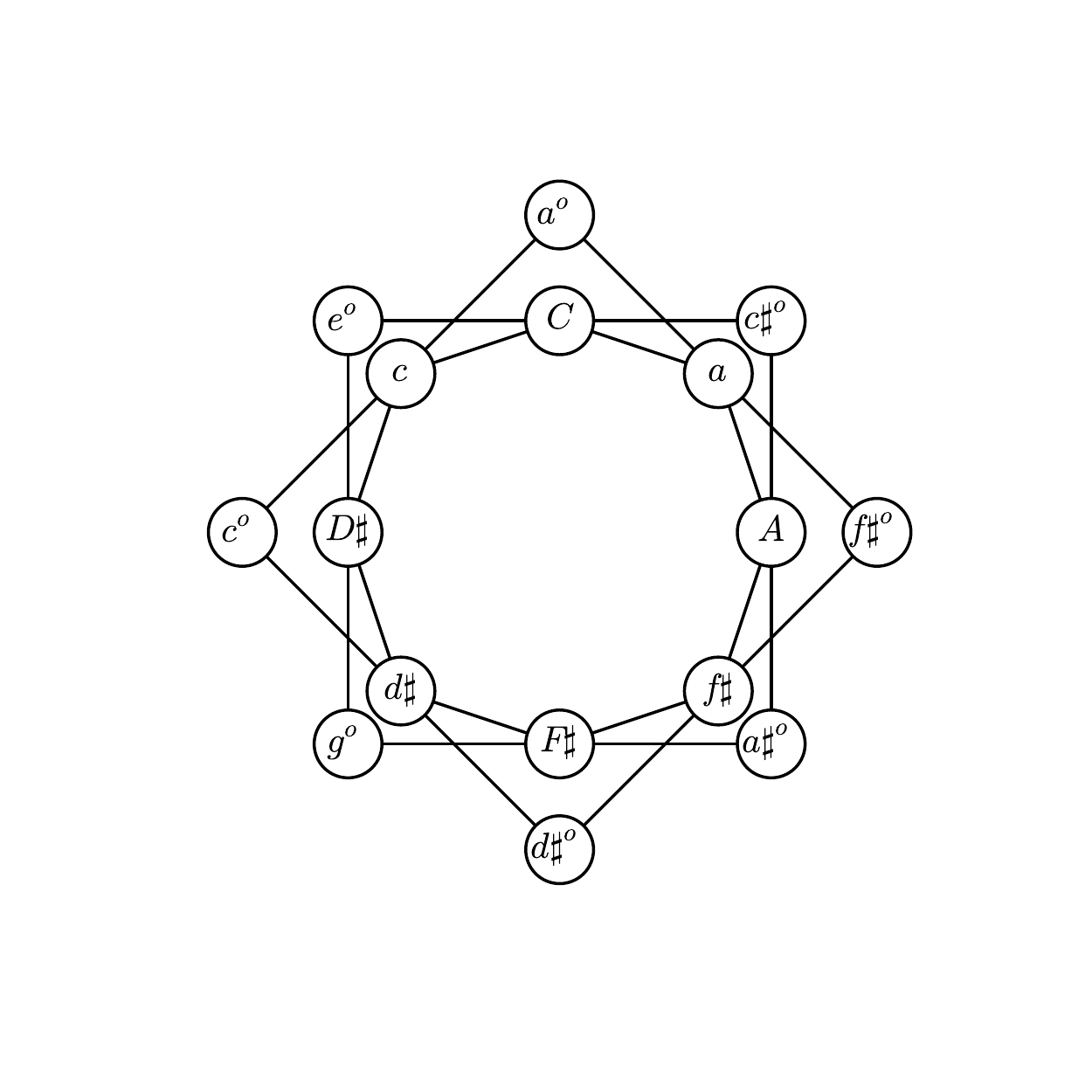}\label{fig:octatonic2}}
\caption{Two views of the voice-leading relations among triads constructable from the octatonic scale, which reveal the symmetric structure of the symmetric diminished scale.}\label{fig:octatonic}
\end{figure}

Since the graphs we are considering are undirected we would be able to move in either direction between a pair of chords, but then we would not be able to use the reverse progression within the same Eulerian trail. In contrast with Hamiltonicity, it is very easy to tell whether a graph is Eulerian. It is a well known, and very useful, fact that a graph is Eulerian if and only if every vertex has even degree. It is semi-Eulerian if precisely two vertices of the graph have odd degree.  Thus the hexatonic graph of figure~\ref{fig:hexatonic} is not Eulerian. The `symmetric diminished scale' $\{0,1,3,4,6,7,9,10\}$, sometimes referred to simply as the `octatonic scale', is constructed by alternating between intervals of a semitone and a whole tone. This scale produces a beautiful graph which reflects the internal symmetry of the scale which its name implies. This graph is Eulerian, and is shown in figure~\ref{fig:octatonic}. However actually finding Eulerian trails, and in particular, evaluating the number of Eulerian trails within a graph is certainly non-trivial for large graphs.



\subsection{Connectedness}

As discussed, only the whole tone scales produce graphs which are not connected. For all other voice-leading graphs, an instructive question is {\em how well-connected} they are. This turns out to be somewhat problematic since most, if not all, of the standard tools for measuring `connectivity' (for example, vertex connectivity, algebraic connectivity \cite{fiedler1989laplacian}, the isoperimetric number) depend, to some extent, on the number of vertices of the graph. This makes comparison between graphs difficult, especially as scales with the same number of notes may produce graphs with different numbers of vertices. For this reason we turn our attention towards more specific measures such as ‘communicability’ and ‘centrality’ which focus on the properties of the individual vertices of a graph.


\section{Complex Networks}

Complex network theory is a modern, vibrant area of study which seeks to analyse graphs which represent, among other things, social networks and computer networks. Such graphs are typically gargantuan in size, and display highly complex, and often dynamic, structures. Because of this, many tools have been developed which aim to concisely distill information about both the individual vertices and the network as a whole \cite{newman2008mathematics}. Though our graphs are small compared to a typical complex network, many of these tools still have specific and appropriate interpretations for parsimonious voice-leading.

\begin{figure}
\centering
\includegraphics[width=0.5\linewidth]{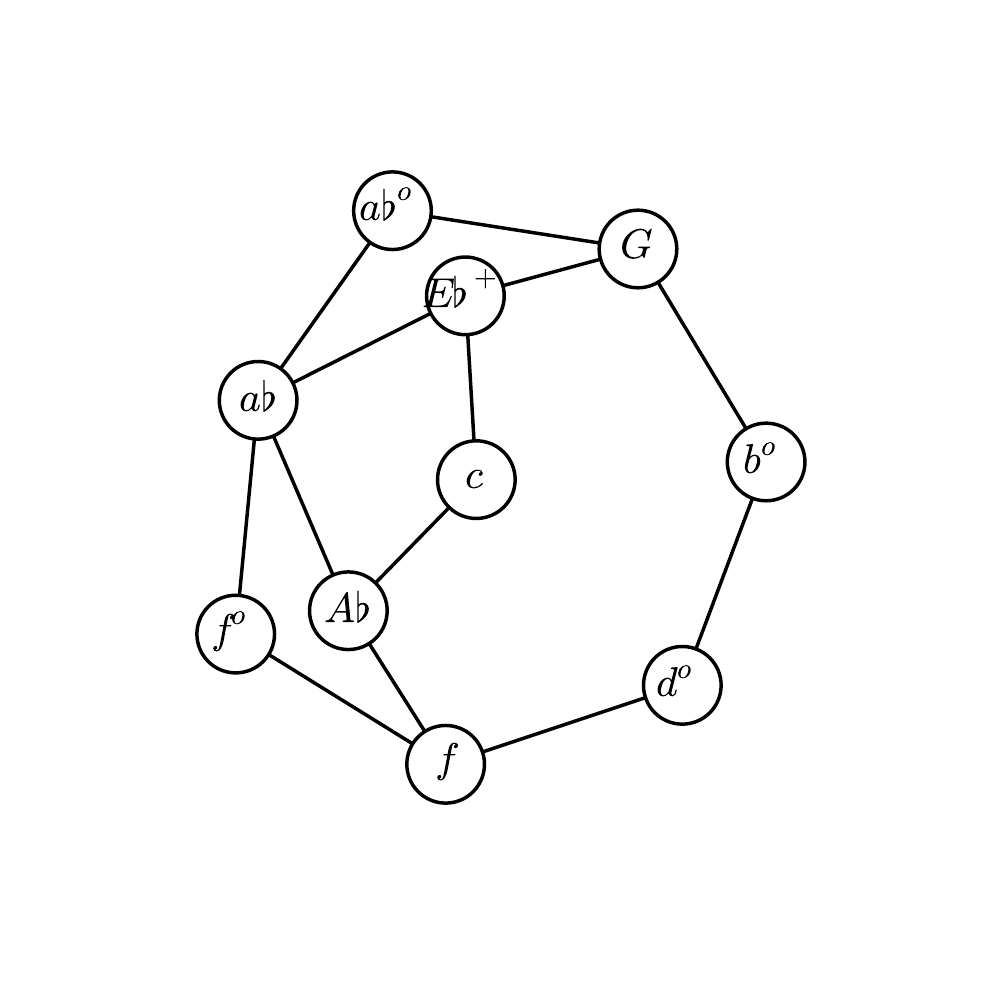}
\caption{Voice-leading triads from scale of $C$ harmonic minor.}\label{fig:Charm_minor_graph}
\end{figure}

Measures of \textit{centrality} determine how `important' a certain vertex is within a network, where a variety of different meanings of `important' have lead to a variety of different meaures. Here we describe several, and discuss their application to voice-leading, using the particular (simple, but instructive) example of triads available in the harmonic minor scale on $C$. The resulting graph for this set of pitches, $\mathcal{P} = \{ 0,2,3,5,7,8,11  \}$, is shown in figure~\ref{fig:Charm_minor_graph}.

The simplest such measure is to consider the degree distribution of the vertices. This gives an absolute value to the amount of harmonic choice at each vertex --- the number of triads that are a single voice-leading step away. For a graph of $n$ vertices, this is commonly normalised by dividing by $(n-1)$, the maximum possible degree at each vertex, to give {\em degree centrality}. Table~\ref{tab:Katz} gives the degree centrality for each of the vertices in figure~\ref{fig:Charm_minor_graph}, which simply indicates what proportion of all available triads are accessible from each triad.

\begin{table}
\begin{center}\caption{A variety of centrality measures for the voice-leading graph of triads from $C$ harmonic minor, shown in figure \ref{fig:Charm_minor_graph}. The values for Katz centrality are computed with $\alpha = 0.35$. This is close to the maximum possible value of $\alpha$, which is given by the reciprocal of the largest eigenvalue of $A$, which is approximately equal to  2.768.}\label{tab:Katz}
\begin{tabular}{|c|c|c|c|c|c|}
\hline
Triads & Degree& Degree  & Closeness & Betweenness &Katz  \\
& & centrality & centrality&centrality&centrality\\
\hline
$c$ minor&2 & 2/9&0.45&0.032&0.276\\
$d$ dim., $b$ dim. &2& 2/9 &0.45&0.111&0.180\\
$E\flat$ aug., $A\flat$ major &3& 3/9& 0.529&0.157&0.380\\
$f$ minor, $G$ major&3 & 3/9 &0.500&0.190&0.305\\
$f$ dim., $a\flat$ dim. &2& 2/9& 0.474&0.051&0.283\\
$a\flat$ minor&4 &4/9& 0.563&0.255 & 0.474\\
\hline
\end{tabular}

\end{center}
\end{table}

The degree of a vertex gives purely local information; it says nothing about how far away other triads are. For a given vertex, one might average the geodesic distance to all other nodes. A smaller value of this average would indicate a greater centrality, so the reciprocal defines the next measure we discuss.

\begin{definition}
The \textbf{closeness centrality} \cite{brandes2007centrality} of a vertex $v_i$ is given by
$$
C(v_i) = \left( \frac{1}{n-1}\sum_{v_j \ne v_i} {d(v_i,v_j)}\right)^{-1}
$$
\end{definition}

Referring to table \ref{tab:Katz}, we can see that closeness centrality distinguishes between, for example, $c$ minor and $f^o$, which have the same degree centrality. On average, fewer chord progressions are necessary to explore the harmonic landscape from $f^o$ than from $c$ minor. 

Given that geodesic distances represent efficient voice-leading, one might ask which vertices appear most often in shortest paths, that is, which triads frequently act as musical stepping-stones between others in voice-leading. This is the quantity measured in the following.
\begin{definition}The \textbf{betweenness centrality} \cite{freeman1977set,brandes2001faster} of a given vertex $v$ is the proportion of shortest paths in the graph which pass through $v$.
\end{definition}

Consider the triads $E\flat^+$ and $f$ in table~\ref{tab:Katz}. $E\flat^+$ is more central, in the sense that it is, on average, a shorter distance from all others than $f$, but $f$ is of more use when constructing the most efficient route between triads. The chord $c$ has the uniquely smallest value of betweenness centrality, which speaks to the fact that we might be more likely to use the C minor chord at the beginning of end of a chord progression, rather than in the middle.

For a final centrality measure, we consider {\em Katz centrality}, which generalises the idea of degree centrality from the local to the global. 
\begin{definition}
Given a network, the \textbf{Katz centrality} \cite{katz1953new} of $v_i$ is given by $$K(v_i)=\sum_{k=1}^\infty \sum_{j=1}^n \alpha^k (A^k)_{ji}$$
where $n$ is the number of vertices in the network, $A$ is the adjacency matrix and $\alpha$ is an attenuation factor which specifies the relative significance to be placed on nearby neighbours within a graph, as opposed to those which are further away. To ensure that the sum converges, $\alpha$ must be smaller than the reciprocal of the largest eigenvalue of $A$.
\end{definition}
Recall that the $ij$-th entry of $A^k$ gives the number of walks of length $k$ from $v_i$ to $v_j$. Thus Katz centrality measures the number of all vertices that can be reached, with the contribution from vertices $k$ steps away diminished by a factor of the $k$th power of $\alpha$. Musically, while degree centrality represents harmonic choice at a given vertex, Katz centrality measures total harmonic choice along all progressions from the initial triad. In terms of social network analysis this effect is often described as `it's not how many friends you have, but how influential your friends are.' This interpretation can again be illustrated by table \ref{tab:Katz}. Triads $f^o$ and $d^o$ have the same degree (and note that $d^o$ has the higher betweenness centrality), but the Katz centrality for $f^o$ is higher, since it has $a\flat^o$, the most influential and central vertex of all, as an immediate neighbour.

As an alternative to looking at individual triads, we might consider how well specific pairs of vertices `communicate' with each other.
\begin{definition}Let $P_{v_iv_j}^{(s)}$ be the number of shortest paths with length $s$ between $v_i$ and $v_j$  and let $W_{v_iv_j}^{(k)}$ be the number of walks between $v_i$ and $v_j$ with length $k>s$. The \textbf{communicability} \cite{estrada2009communicability} between the vertices $v_i$ and $v_j$ is then
$$C_{v_iv_j}=\frac{1}{s!}P_{v_iv_j}+\sum_{k>s}\frac{1}{k!}W_{v_iv_j}^{(k)}.$$
\end{definition}

This measure considers all routes between $v_i$ and $v_j$, with shorter routes being weighted more heavily than longer ones. The communicability of a pair of vertices then captures the idea of how easy it is to move between these two vertices. Using the $k$th power of the adjacency matrix, this expression can be rewritten as $$C_{v_iv_j}=\sum_{k=0}^{\infty}\frac{(A^k)_{v_iv_j}}{k!}=(e^A)_{ij}.$$

Again we can consider the $C$ harmonic minor example. Although $f$ and $a\flat$ are both a distance of two away from $c$, we have $C_{cf}=0.998$ and $C_{c a\flat}=1.859$. It is also interesting to consider the values when each vertex is paired with itself. In this case $a\flat$ has the highest value, at 4.500, with $b^o$ and $d^o$ having the lowest values at 2.343. This agrees with the intuition, since the more isolated nature of these two vertices means that a smaller proportion of walks, or chord progressions, will use these vertices. 

\section{Conclusion and Extensions}

We have discussed the application of a variety of concepts and tools from graph theory and complex network theory to efficient voice-leading within a given harmonic structure. The notion of a self-centred graph demonstrated that in one particular (classical graph-theoretic) sense, typical musical scales are all equally well-distributed, in that no triad is any more or less isolated than any other. However, the relative centrality of vertices can be distinguished using modern ideas from complex networks, each having a specific interpretation to harmonic progressions. This field of study is an evolving area, with a  potentially increasing number of new tools which could be incorporated into this framework. 

There are natural musical extensions to this approach, for example exploring the connectivity of voice-leading in $n$-tone equal temperament, or unequal temperament. There are no mathematical issues with considering such graphs, indeed, many of the measures discussed here are designed for larger graphs. Weighted graphs could be employed to represent notions of voice-leading size, and similarly directed graphs could model situations in which chord progressions are allowed in only one direction. Many of the graph theoretic tools we have employed have direct analogues for weighted and directed graphs so the analysis should not prove too difficult in these cases. A further idea of interest would be to consider dynamic graphs, changing with time to represent the idea of changing key signatures or tonal centres throughout a piece. In terms of musical analysis and composition, it is hoped that these ideas offer a new viewpoint from which to study parsimonious voice-leading, and even give inspiration to composers who seek to use voice-leading in some efficient manner.


\end{document}